\documentclass[12pt,fleqn]{article}
\usepackage{amssymb,amsmath,graphicx,cite}
\usepackage[cp1251]{inputenc}

\begin{document}

\begin{center}
\LARGE \bf Antimonotonous quadratic forms \\
and partially ordered sets
\end{center}

\begin{center}
\Large\bfseries L.A. Nazarova~$^*$, A.V. Roiter~$^{*}$, M.N. Smirnova~$^{*}$
\end{center}

{ \footnotesize

$^{*}$~Institute of Mathematics of National Academy of Sciences of Ukraine,\\
Tereshchenkovska str., 3, Kiev, Ukraine, ind. 01601\\
E-mail: roiter@imath.kiev.ua}

\bigskip

Let $P$ be a bounded set in $n$-dimensional space
${\mathbb R}_n$, $f(x_1,\ldots, x_n)=f(x)$ ($x\in \mathbb{R}_n$) be a
continuous function.
By the second Weierstrass theorem (\cite{13}, s.~163) $\inf
\{f(\overline{P})\}$ $(=\inf\limits_{\overline P}f(x))$ is reached. We will
call a function $f$ {\it $P$-faithful}, if $\inf
\{f(\overline{P})\}$
is not reached on $\overline{P}\setminus P$ and
$\inf\{f(\overline{P})\}>0$ (i.~e. $f$ is positive on $\overline P$).
Remark that if $n=1$, $P=(a,b)$, then it is necessary for $P$-faithfulness
that a function is not monotonous.

Further we will suppose that $P=P_n=\{(x_1,\ldots,x_n)\,|\,0<
x_i\leq 1, \,i=\overline{1,n},\, x_1+\cdots+x_n=1\}$. If $n>1$ then
$x_i<1$, $i=\overline{1,n}$. Then $\overline{P}=\{(x_1,\ldots,x_n) | 0\leq
x_i\leq 1,\; i=\overline{1,n},\; x_1+\ldots +x_n=1\}$. In this case the
$P$-faithfullness $f$ is essentially connected with behavior of function on
the hyperplane $H_n=\{(x_1,\ldots,x_n)\,|\,
x_1+\cdots+x_n=0\}$.

For differentiable function $f$
put $C^-(f)=\{h\in H_n\setminus \{0\}\,|\, \frac{\partial f}{\partial
x_i}(h)\leq 0,\, i=\overline{1,n}\}$, $C^+(f)=\{h\in H_n\setminus
\{0\}\,|\, \frac{\partial f}{\partial x_i}(h)\geq 0,\, i=\overline{1,n}\}$,
$C(f)=C^+(f)\cup C^-(f)$. A function $f$ is {\it antimonotonous} if
$C(f)=\varnothing$. If $n=1$ then $P_1=(1)$, $H_1\setminus
\{0\}=\varnothing$ and any function is antimonotonous.

In s.~3 (Proposition 1) we prove that any $P$-faithful quadratic form
is antimonotonous and, therefore, in this case antimonotonousness is a
generalization of $P$-faithfulness.

\medskip

{\bf Example 1.} A linear function $f=\sum\limits_{i=1}^n a_i x_i$ is
antimonotonous only if all $a_i=0$, i.~e.
$f=0$.
The quadratic forms $x_1^2+x_2^2$, $x_1^2+x_2^2+x_1x_2$ are antimonotonous,
but the forms $x_1^2-x_2^2$, $x_1^2+x_2^2+x_3^3+x_1x_2+x_1x_3$ are not.

\medskip

The problem of effective criterion of antimonotonousness even for
quadratic forms is probably difficult.

In this article we solve this problem for quadratic form $f_S$,
attached to (finite) partially ordered set (poset)
$S=\{s_1,\ldots,s_n\}$: $f_S(x_1,\ldots,x_n)=\sum\limits_{s_i\leq s_j}x_i x_j$ \cite{1}
under the additional condition of positive semidefiniteness
form $f_S$ (i.~e. $f_S(x)\geq 0$).
Posets with antimonotonous form
generalize $P$-faithful posets, defined in \cite{2}
and studied in \cite{2,3,4,5,6}
and (it follows from this work) coincide with them not only for positive
definite forms, but for positive semidefinite ones.

A direct construction of a vector contained in $C(f_S)$ permits essentially
to simplify the proof of criterion of $P$-faithfulness \cite{2,4,5,6}
avoiding consideration of many different cases.

We also take out evident formula for calculating of
$\inf\{f_S(\overline{P})\}$ for $P$-faithful $S$, on base of this formula
we give simple proofs of the criterions of finite representativity \cite{7}
(see also \cite{8} and tameness \cite{9}, (see also \cite{10}) of partially
ordered sets.

{\bf 1.} In s.~1 $f$ is a
differentiable function, defined on ${\mathbb R}_n$.
We call vectors the elements of ${\mathbb R}_n$.

Let ${\mathbb
R}^+_n=\{x\in{\mathbb R}_n\,|\,
x_i>0,\,i=\overline{1,n}\}$, ${\overline{\mathbb
R}}^+_n=\{x=(x_1,\ldots,x_n)\in {\mathbb R}_n\,|\, 0\leq x_i,
i=\overline{1,n}; x\not=0\}$, $P_n=\overline{P}_n\cap {\mathbb R}^+_n$
($\mathbb{R}^+=\mathbb{R}_1^+$).

If $f_1$, $f_2$ are defined respectively on ${\mathbb R}_m$
and on ${\mathbb R}_n$
we put $(f_1\oplus f_2)(x_1,\ldots,x_m$,
$x_{m+1},\ldots,x_{n+m})=f_1(x_1,\ldots,x_{m})+
f_2(x_{m+1},\ldots,x_{m+n})$.

We call a twice differentiable function $f$ {\it concave},
if the following conditions hold:

a) $\frac{\partial f}{\partial x_i}(0)=0,\; i=\overline{1,n}$,

b) ${\partial^2f\over \partial x_i \partial x_j} \geq 0$,
$i,j=\overline{1,n}$,

\noindent
{\it $q$-conceve}, $q\in\mathbb{R}^+$, if a), b) and

c) ${\partial^2f\over \partial x_i^2} \geq q$,
$i=\overline{1,n}$ hold.

\noindent
A quadratic form $f_S$ attached to poset $S$ is, in particular, 2-concave.

\medskip

{\bf Remark 1.} By Lagrange theorem
b) implies $I$,  c) implies II$_{q}$.

\medskip

I) $\frac{\partial f}{\partial x_i}(x_1,\ldots,x_{j-1},x_j+d,x_{j+1},\ldots,x_n)\geq
\frac{\partial f}{\partial x_i}(x_1,\ldots,x_n)$, $(i, j\in 1,\ldots,n)$.

II$_q$) $\frac{\partial f}{\partial x_i}(x_1,\ldots,x_{i-1}, x_{i}+d,x_{i+1},\ldots,x_n)\geq
\frac{\partial f}{\partial x_i}(x_1,\ldots,x_n)+qd$, $(i=\overline{1,n})$;

Put $\widehat C^-(f)=\Big\{x\in {\mathbb R}_n\setminus
\{0\}\,|\,\sum\limits_{i=1}^n x_i\geq 0,\ {\partial f\over\partial
x_i}(x)\leq 0,\ i=\overline{1,n}\Big\}$; $\widehat C^+(f)=\Big\{x\in
{\mathbb R}_n\setminus \{0\}\,|\,\sum\limits_{i=1}^n x_i\leq 0,\ {\partial
f\over\partial x_i}(x)\geq 0,\ i=\overline{1,n}\Big\}$.

\medskip

{\bf Lemma 1.} {\it If $f$ is a concave function then it is
antimonotonous if and only if $\widehat  C^+(f)\cup \widehat
C^-(f)=\varnothing$.}

\medskip

Let $x\in \widehat C^-(f)$ (the case $x\in \widehat C^+(f)$ is
analogous), $\sum\limits_{i=1}^d x_i=d\in {\mathbb R}^+$. Then
$\{x_1-d,x_2,\dots,x_n\}\in C(f)$ (according to I) if only $x\neq
(d,0,\dots,0)$. But at the last case $(d,-d,0,\dots,0)\in C(f)$. If
$y\in C(f)$ then it is clearly that $y\in \widehat C^+(f)\cup \widehat
C^-(f)$. \hfill \rule{2mm}{2mm}

\medskip

{\bf Lemma 2.} {\it If $f_1$ and $f_2$ are concave then the function
$f_1\oplus f_2$ is antimonotonous if and only if $f_1$ and $f_2$ are
antimonotonous.}

\medskip

We shall prove that $C(f_1\oplus f_2)\neq\varnothing$ if and only if
either $C(f_1)\neq\varnothing$ or $C(f_2)\neq\varnothing$. If
$(x_1,\dots,x_{n_1},y_1\dots,y_{n_2})\in C(f_1\oplus f_2)$ then either
$(x_1,\dots,x_{n_1})\in \widehat C^+(f_1)\cup \widehat C^-(f_1)$ or
$(y_1,\dots,y_{n_2})\in \widehat C^+(f_2)\cup \widehat C^-(f_2)$ and by
Lemma~1 in the first case $C(f_1)\neq\varnothing$ and in the second case
$C(f_2)\neq \varnothing$.

If $x=(x_1,\dots,x_{n_1})\in C(f_1)$ then
$(x_1,\dots,x_{n_1},\underbrace{0,\dots,0}\limits_{n_2})\in
C(f_1\oplus f_2)$ using $a)$; if $(y_1,\dots,y_{n_2})\in
C(f_2)$ then
$(\underbrace{0,\dots,0}\limits_{n_1},y_1,\dots,y_{n_2})\in
C(f_1\oplus f_2)$.
\hfill\rule{2mm}{2mm}

\medskip

We say that a nonzero vector $d\in \mathbb{Z}_n$ is {\it $m$-Dynkin})
$(1\leq m\leq n)$
for $q$-concave function $f$, if
1) $0\leq \frac{\partial f}{\partial x_m}(d)\leq q$;
2) $\frac{\partial f}{\partial x_j}(d)=0$ if $j\not= m$,
$j=\overline{1,n}$.

We call a function $f$ {\it $m$-isolated}, if
${\partial f\over\partial x_k}(s_m)=0$ for $1\leq
k\leq n$, $k\neq m$,
$s_m=(\underbrace{0,\dots,0}\limits_{m-1},1,0,\dots,0).$

\medskip

{\bf Lemma 3.} {\it  Let $f$ be $q$-concave not $m$-isolated function.
If there exists $m$-Dynkin vector then $C(f)\not=\varnothing$.}

\medskip

Let $\sum\limits_{i=1}^n d_i=\overline d$.
If $\overline d\leq 0$ then $d\in \widehat C^+(f)$
and $C(f)\not=\varnothing$
by Lemma~1. Let $\overline d>0$. Put $u_j=d_j$ if $j\not =m$ and
$u_m=d_m-\overline d$, and we prove that $u=(u_1,\ldots,u_n)\in C(f)$. It
is clear that $u\in H_n$.
$\frac{\partial f}{\partial x_j}(u)\leq 0$ if $j\not=m$
according to I and 2, and $\frac{\partial f}{\partial x_m}(u)\leq 0$
according to II$_q$ and 1).

Now we shall prove that $u\neq 0$. If $u=0$ then $d=\lambda s_m$,
$\lambda\neq 0$ (since $d\neq 0$). But then non $m$-isolateness
$f$ implies that ${\partial f\over\partial
x_k}(d)\neq 0$, $k\neq m$.\hfill\rule{2mm}{2mm}

\medskip

{\bf Example 2.} Let $S=\{s_1,s_2,s_3,s_4,s_5\,|\,
s_1<s_i,\,i=\overline{2,5}\}$,
$f_S=\sum\limits_{i=1}^5 x_i^2+x_1\sum\limits_{j=2}^5 x_j$,
$d=(-2,1,1,1,1)$ --- $i$-Dynkin vector
for $f_S$ $(i=\overline{1,5})$.
Vectors $(-2,1,1,1,-1)$,
$(-2,1,1,-1,1)$, $(-2,1,-1,1,1)$, $(-2,-1,1,1,1)$ and $(-4,1,1,1,1)$
belong to $C(f_S)$.

\medskip

Consider $P$-faithfulness.
Put ${\rm St}(f)=\{a\in{\mathbb R}^+_n\,|\,
\frac{\partial f}{\partial x_i}(a)=\frac{\partial f}{\partial x_j}(a)$,
$i,j=\overline{1,n}\}$;
${\rm St}^+(f)=\{a\in{\rm St}\,|\,
\frac{\partial f}{\partial x_i}(a)>0\}$.

Vector $u\in P_n$ is {\it $P$-faithful} for function $f$,
if $f(u)>0$ and $w\in \overline{P}_n$
implies $f(u)\leq f(w)$,
moreover if $w\not\in P_n$,
then $f(u)<f(w)$.

Denote $\widetilde{{\rm St}}(f)$ the set of
$P$-fathful vectors for $f$.
$P$-faithfulness of $f$ is equivalent to $\widetilde{{\rm St}}(f)\not=0$.

\medskip

{\bf Lemma 4.} {\it $\widetilde{{\rm St}}(f)\subseteq
{\rm St}(f)$ for any $f$.}

\medskip

Let $n>1$. We express $x_n=1-\sum\limits_{i=1}^{n-1}x_i$, and obtain the
function $\widehat{f}(x_1,\ldots,x_{n-1}) = f(x_1,\ldots, x_{n-1},
1-\sum\limits_{i=1}^{n-1}x_i)$. If $u=(u_1,\ldots,u_n)$ is $P$-faithful
vector for $f$ then $\widehat{f}(\widehat{u})$ is minimum when
$\widehat{u}=(u_1,\ldots,u_{n-1})$. $\frac{\partial \hat f}{\partial
x_i}=\frac{\partial f}{\partial x_i}+ \frac{\partial f}{\partial x_n}\cdot
\frac{\partial x_n}{\partial x_i}$; $x_n=1-\sum\limits_{i=1}^{n-1}x_i$;
$\frac{\partial x_n}{\partial x_i}=-1$ $(i=\overline{1,n-1})$. Therefore
$\frac{\partial \hat f}{\partial x_i}=\frac{\partial f}{\partial x_i}-
\frac{\partial f}{\partial x_n}=0$. \hfill \rule{2mm}{2mm}

\medskip

Let $f$ be a homogeneous function of degree $k$ (i.~e.
$f(\lambda x_1,\ldots,\lambda x_n)=\lambda^k f(x_1,\ldots,x_n)$).
For homogeneous $f$, if $k\not=1$ and $\inf\{f(P)\}>0$ put $P(f)=\inf
\{f(P)\}^{\frac{1}{1-k}}$. In particular, for $k=2$, $P(f)=\inf
\{f(\overline{P})\}^{-1}$.

\medskip

{\bf Lemma 5.} {\it Let $f_1(x_1,\ldots,x_{n_1})$, $f_2(x_{n_1+1},\ldots, x_{n_2})$
be two homogeneous functions of degree $k$,
$n_1+n_2=n$, $\inf\{f_j(P_{n_j})\}>0$, $j=1,2$.
Then $P(f_1\oplus f_2)=P(f_1)+P(f_2)$.}

\medskip

Values of homogeneous function $f$ on ${\overline{\mathbb
R}}^+_n$ are defined by its values on $\overline{P}_n$, namely for $y\in
{\overline{\mathbb R}}^+_n$, $f(y)=\lambda^kf(u)$ where $u\in
\overline{P}_n$, $\lambda=\sum\limits_{i=1}^ny_i$, $u=\lambda^{-1}y$.
Therefore $\inf\{(f_1\oplus
f_2)(\overline{P}_n)\}=\inf\limits_{0\leq\lambda\leq 1}\bigg[\lambda^k
\inf\{f_1(\overline{P}_{n_1})\}+(1-\lambda)^k\inf\{f_2(\overline{P}_{n_2})\}
\ b i g g ] $ . Put $\inf\{f_1(\overline{P}_{n_1})\}=a$,
$\inf\{f_2(\overline{P}_{n_2})\}=b$.

Consider function $\Phi_{ab}(\lambda)=a\lambda^k+b(1-\lambda)^k$,
$a>0$, $b>0$ and find $\inf\limits_{0\leq\lambda\leq
1}\Phi_{ab}(\lambda)$.
The derivative of $\Phi_{ab}(\lambda)$ with respect to $\lambda$ (consider
$u$ and $v$ to be constants) is
$(\Phi_{ab}(\lambda))'_\lambda=ka\lambda^{k-1}-kb(1-\lambda)^{k-1} $ . Let
$\overline\lambda$ be a positive root of equation
$(\Phi_{ab}(\lambda))'_\lambda =0$. Substitute $\overline\lambda$ and
obtain $a\overline\lambda^{k-1}=b(1-\overline\lambda)^{k-1}$. Raise the
both parts of this equality to ${1\over{k-1}}$ power and obtain
$a^{{1\over{k-1}}}\overline\lambda= b^{1\over{k-1}}(1-\overline\lambda)$,
$\overline\lambda = {b^{1\over{k-1}}\over
a^{1\over{k-1}}+b^{1\over{k-1}}}$. Thus, $\inf\limits_{0\leq\lambda\leq 1}
\Phi_{ab}(\lambda)=\min\{\Phi_{ab}(0),\Phi_{ab}(1),
\Phi_{ab}(\overline\lambda)\}=\min\{a,b,\Phi_{ab}(\overline\lambda)\}$.

We will prove that $\Phi_{ab}(\overline\lambda)<\Phi_{ab}(0)=b$.
Really,
$a\overline\lambda^{k-1}=b(1-\overline\lambda)^{k-1}$.
$\Phi_{ab}(\overline\lambda)=a\overline\lambda^k+b(1-\overline\lambda)^k=
b(1-\overline\lambda)^{k-1}\overline\lambda+b(1-\overline\lambda)^k=
b(1-\overline\lambda)^{k-1}$.
$\Phi_{ab}(\overline\lambda)<b$, since $a>0$, $b>0$, $0<\overline\lambda<1$
Analogously
$\Phi_{ab}(\overline\lambda)=a\overline{\lambda}^{k-1}
<\Phi_{ab}(1)=a$. Therefore
$\inf\{(f_1\oplus f_2)(\overline{P}_n)\}=\Phi_{ab}(\overline\lambda)=
{ab\over\left(a^{1\over k-1}+b^{1\over k-1}\right)^{k-1}}$.

Let us return to $P(f_1\oplus f_2)$, $P(f_1)$, $P(f_2)$. We have
$P(f_1)=a^{1\over{1-k}}$, $P(f_2)=b^{1\over 1-k}$,
$P(f_1\oplus f_2)=\Bigg({ab\over\left(a^{1\over k-1}+b^{1\over
k-1}\right)^{k-1}}\Bigg)^{1\over 1-k}=b^{1\over 1-k}+a^{1\over
1-k}$. \hfill\rule{2mm}{2mm}

\medskip

{\bf Corollary 1.} {\it Under the conditions of the lemma $f_1\oplus f_2$
is $P$-faithful if and only if $f_1$ and $f_2$ are $P$-faithful.}

\medskip

{\bf 2.} Further $f=\sum\limits_{i,j=1}^n a_{ij} x_i x_j$
$(a_{ij}=a_{ji})$ is a quadratic form over field
${\mathbb R}$;
$A=(a_{ij})$
is a symmetric matrix of quadratic form $f$.
$\frac{\partial^2 f}{\partial x_i \partial x_j} = 2a_{ij}$. $f$ is a
$2$-concave if $a_{ii}\in{\mathbb N}$,
$a_{ij}+a_{ji}\in{\mathbb N}_0$ $(i,j=\overline{1,n})$. For
$x=(x_1,\ldots,x_n)\in{\mathbb R}_n$ (for fixed $f$) put
\[
x'_i=\frac{\partial f}{\partial
x_i}(x_1,\ldots,x_n)=2\sum_{j=1}^na_{ji}x_j, \quad
x'=(x'_1,\ldots,x'_n)=2xA.
\]

We will use the following identity, which is
easy to check: \begin{gather} f(u+v)=f(u)+f(v)+\sum_{i=1}^n u'_i v_i, \quad
u,v\in {\mathbb R}_n,\quad\mbox{that implies}\label{0*}\\ \sum_{i=1}^n u'_i
v_i=\sum_{i=1}^n v'_i u_i, \quad f(u+\varepsilon
v)=f(u)+\varepsilon^2f(v)+\varepsilon \sum_{i=1}^n u_iv'_i, \quad
\varepsilon\in{\mathbb R}.\label{*} \end{gather}

Put in \eqref{0*} $u=v$, and obtain (see.~\cite{13}, s.~178):

\begin{equation}\label{*1}
f(u)=\frac 12 \sum_{i=1}^n u_i u'_i.
\end{equation}

By means of (3) we can reformulate Lemma~4 by the following way.

\medskip

{\bf Lemma 4$'$.} {\it $\widetilde{{\rm
St}}(f)\subset {\rm St}^+(f)$ for quadratic form $f$.}

\medskip

Put
$\widetilde{C}(f)=\{(v_1,\ldots,v_n)\in C(f)\,|\,
(v'_1,\ldots,v'_n)\not=0\}$. Since $\frac{\partial f}{\partial
x_i}(-x)=-\frac{\partial f}{\partial x_i} (x)$ then if $C(f)\neq
\varnothing$, then $C^-(f)\neq \varnothing$ and $C^+(f)\neq \varnothing$.
Therefore, choosing a vector $v\in C(f)\neq\varnothing$, further we will
suppose that $v\in C^-(f)$.

\medskip

{\bf Proposition 1.} {\it For any quadratic form  $f$
$1)$ at least one of
the sets ${\rm St}(f)$ and $\widetilde{C}(f)$ is empty
$2)$ at least one of the sets $C(f)$ and $\widetilde{\rm St}(f)$ is empty.}

\medskip

1) Let $u\in {\rm St}(f)$, $v\in \widetilde{C}(f)$. Then
$\sum\limits_{i=1}^nu'_iv_i=u'_1\sum\limits_{i=1}^nv_i=0$. In the other
hand, if $v'_j<0$ then $\sum\limits_{i=1}^nu_iv'_i<0$ (since $v'_i\leq 0$,
$u_i>0$, $i=\overline{1,n}$), that contradicts~\eqref{*}.

2) Let $u\in\widetilde{\rm St}(f)$, $v\in C(f)$. If
$v\in\widetilde{C}(f)$ then 1) implies the statement. Let $v\in
C(f)\setminus\widetilde{C}(f)$, i.~e. $v'_i=0$ for
$i=\overline{1,n}$.
Then $f(v)=0$ (\eqref{*1}) and so
$f(u+\varepsilon v)=f(u)$
for any $\varepsilon$.
Put $|\varepsilon|=\min\limits_i\frac{u_i}{|v_i|}$.
The sign of $\varepsilon$
is opposite sign of one of those $v_i$,
for which the minimum is reached.
Then $u+\varepsilon v\in\overline{P}_n\setminus P_n$,
that contradicts $P$-faithfulness of $u$.~\hfill~\rule{2mm}{2mm}
\medskip

{\bf Corollary 2.} {\it $P$-faithful quadratic form is antimonotonous.}
\medskip

{\bf Example 3.} Let
$f=\sum\limits_{i=1}^4x_i^2+(x_1+x_2)(x_3+x_4)$. Then
\[
A= \left(\begin{array}{cccc}
1&0&{1\over 2}&{1\over 2}\\
0&1&{1\over 2}&{1\over 2}\\
{1\over 2}&{1\over 2}&1&0\\
{1\over 2}&{1\over 2}&0&1
\end{array}\right),\quad St(f)\ni (1,1,1,1),\quad
C(f)\ni (1,1,-1,-1).
\]
Proposition~1 implies $\widetilde{{\rm St}}(f)=\varnothing$,
$\widetilde{C}(f)=\varnothing$.

In this example $|A|=0$.

\medskip

{\bf Proposition 2.} {\it If $|A|\neq 0$ then
one of the sets $C(f)$, ${\rm St}(f)$ is not empty, but the other one is
empty.}

\medskip

At first suppose that $\varnothing\neq C(f)\ni v$ and
$\varnothing\neq {\rm St}(f)\ni u$. If $v\in\widetilde{C}(f)$ then
${\rm St}(f)=\varnothing$ by Proposition~1. If $v\in
C(f)\setminus\widetilde{C}(f)$ then $v'=0$, and $vA={1\over 2}v'=0$.
Hereof $v=0$ that contradicts $v\in H_n\setminus \{0\}$ (the
definition of $C(f)$).

Now we prove that either ${\rm St}(f)\neq\varnothing$, or
$C(f)\neq\varnothing$. Let $e_n=(1,\dots,1)\in \mathbb{R}_n$,
$y=e_nA^{-1}$, $yA=e_n$. If $y\in \mathbb{R}_n^+$ or $-y\in \mathbb{R}_n^+$,
then $y\in {\rm St}(f)$. If $\{y,-y\}\cap \mathbb{R}_n^+=\varnothing$ then
either for certain $k$ $y_k=0$, or for certain $s$ and $t$ $y_s<0$,
$y_t>0$. It is easy to see that in both cases there exists
$w\in\overline{\mathbb{R}}_n^+$ such that
$wy^T
\Big(=\sum\limits_{i=1}^nw_iy_i\Big)=0$ (at the first case we can put
$w_k>0$, $w_i=0$ for $i\neq k$, at the second case $w_s=y_t$, $w_t=-y_s$,
$w_i=0$ for $i\not\in \{s,t\}$, $i=\overline{1,n}$). We prove that
$v=-wA^{-1}\in C(f)$. $-v'=wA^{-1}A=w\in \overline{\mathbb R}^+_n$, and so
$v_i'\leq 0$. $v\neq 0$, since $w\neq 0$ and $|A|\neq 0$. It remains to
prove that $v\in H_n$, that is equivalent to $v e_n^T=0$. $v
e_n^T=-wA^{-1}e_n^T$; $y^T=(A^{-1})^Te_n^T=A^{-1}e_n^T$ (since $A^T=A$).
Therefore $-wA^{-1}e_n^T=-wy^T=0$. \hfill \rule{2mm}{2mm}

\medskip

{\bf Proposition 3} (see~\cite{6}, p.~II, remark to theorem~1). {\it
$1)$ If $\widetilde{{\rm St}}(f)\not=\varnothing$,
then $f$ is positive definite.
$2)$ If $f$
is positive definite then $\widetilde{{\rm St}}(f)={\rm St}(f) \cap
P_n$ (and therefore $\widetilde{{\rm St}}(f)=\varnothing$ if and
only if ${\rm St}(f)=\varnothing$).}

\medskip

1) We suppose contrary: $f(v)\leq 0$ $(v\not=0)$,
$u\in \widetilde{{\rm St}}(f)$.
a) At first we suppose that $v\in H_n$, i.~e. $\sum\limits_{i=1}^n v_i=0$
and $f(v)<0$.
Then $f(u+\varepsilon v)=f(u)+\varepsilon^2 f(v) +\varepsilon
\sum\limits_{i=1}^n u'_i v_i$. By Lemma~4 $u\in{\rm St}(f)$ and, therefore,
$\varepsilon\sum\limits_{i=1}^n u'_i v_i=0$, i.~e. $f(u+\varepsilon
v)=f(u)+\varepsilon^2 f(v)$. Since $f(v)<0$ then $f(u+\varepsilon
v)<f(u)$, that contradicts $P$-faithfulness of $u$.

b) Now let $v\in H_n$, $f(v)=0$.
Then $f(u+\varepsilon v)=f(u)$
for any $\varepsilon$.
Put $\varepsilon=\min\limits_i\frac{u_i}{|v_i|}$.
The sign of $\varepsilon$
is opposite to the sign of one of those $v_i$,
for which this minimum is reached.
Then $u+\varepsilon v\in \overline{P}_n\setminus P_n$,
again we have the contradiction with $P$-faithfulness of $u$.

c) Let finally $\sum\limits_{i=1}^n v_i\not=0$.
We can admit $\sum\limits_{i=1}^n v_i=1$.
Put $w=u-v$. \eqref{*} for $\varepsilon=-1$
and Lemma~4 imply
$f(w)=f(u)+f(v)-u'$, $u'=u'_i$, $i=\overline{1,n}$.
\eqref{*1} implies $f(u)=\frac{u'_1}{2}$ and, therefore,
$f(w)=f(v)-\frac{u'_1}{2}$.
By Lemma~4$'$ $u'>0$, $f(w)<0$; $w\in H_n$.
Hereat $w\not=0$, since if $w=0$ then $u=v$, according to \eqref{*1}, but
$f(v) \leq 0$. Thus, we reduced c) to a).

2) Let $u\in{\rm St}(f)\cap P_n$, $v\in \overline{P}_n$,
$v\not=u$. $u\neq 0$ since $u\in {\rm St}(f)$ so $f(u)>0$. We will prove
that $f(u)<f(v)$. $0<f(u-v)\,\overset{\mbox{\footnotesize formula (2)}}{=}\,
f(u)+f(v)-\sum\limits_{i=1}^n u'_i v_i \,\overset{\mbox{\footnotesize
formula (3)}}{=}\, \frac{u'}{2}+f(v) -u'= f(v)-\frac{u'}{2}=f(v)-f(u)$,
i.~e. $f(v)>f(u)$.\hfill \rule{2mm}{2mm}

{\bf 3.} Further we will consider
2-concave form $f_S$ for poset $S=\{s_1,\ldots,s_n\}$,
$f_S=\sum\limits_{s_i\leq s_j} x_ix_j$.
Put $C(S)=C(f_S)$.
$S$ is {\it antimonotonous} if $f_S$ is antimonotonous.

A poset $S$ is {\it $P$-faithful}, if $\widetilde{{\rm
St}}(f_S)\not=\varnothing$. (It is equivalent to the definition of
$P$-faithfulness of the poset from \cite{2}) In this case
Proposition~1 implies $C(S)=\varnothing$. Remark that
$\inf\{f_S(\overline{P})\}>0$ since $a_{ij}\geq 0$, $i,j=\overline{1,n}$,
$A\not= (0)$.

Hasse quiver (orgraph) $Q(S)$ of poset $S$
is a quiver, whose vertices are elements of $S$ and two vertices are
connected by an arrow $s_i\to s_j$ if $s_i<s_j$ and there is no $s_k\in S$
such that $s_i<s_k<s_j$. Drawing lines (edges) instead of arrows, we obtain
(nonoriented) graph Hasse $\Gamma(S)$ of partially ordered set $S$. Finite
poset $S$ is usually depicted by a diagram, i.~e. by graph $\Gamma(S)$
assuming that lesser element is drawn below than greater.

We denote elements of poset $S$ and corresponding elements of $Q(S)$ and
$\Gamma(S)$ by the same symbol.

Put ${\rm St}(S)={\rm St}(f_S)$, $\widetilde{\rm St}(S)=\widetilde{\rm
St}(f_S)$. Let $\mathcal{T}_S(x_1,\ldots,x_n)=\sum\limits_{i=1}^n
x_i^2-\sum\limits_{s_i-s_j} x_ix_j$ be a quadratic Tits form of graph
$\Gamma(S)$ (the second sum is taken by all edges of graph $\Gamma(S)$).
We denote matrix of the form $\mathcal{T}_S$ either $\mathcal{A}$ or
$\mathcal{A}(S)$.

Let $Q$ be a quiver without loops and parallel (i.~e. having the same
origin and the same terminus) pathes with vertices $s_1,\ldots,s_n$.
$\widetilde{Q}$ is a matrix, in which $\widetilde Q_{ij}=1$, if there is an
arrow from $s_i$ to $s_j$, $\widetilde Q_{ij}=0$ in the opposite case
$(i,j=\overline{1,n})$. Then $({\widetilde{Q}}^t)_{ij}$ is equal to
the number of pathes of length $t$ from $s_i$ to $s_j$,
${\widetilde{Q}^n}=0$. If, moreover, $Q=Q(S)$ then
$A=\frac{1}{2}[(E+\widetilde{Q}+\cdots+\widetilde{Q}^{n-1})+(E+\widetilde{Q}
+ \cdots+ \widetilde{Q}^{n-1})^T]$ ($A$ is the matrix of $f_S$). It is easy
to see that
$(E+\widetilde{Q}+\cdots+\widetilde{Q}^{n-1})=(E-\widetilde{Q})^{-1}$. Put
$E-\widetilde{Q}=\hat Q$, $|\hat Q|=1$. $\mathcal{A}=\frac 12 (\hat Q+\hat
Q^T)$. $A=\frac{1}{2}(\widehat{Q}^{-1}+(\widehat{Q}^{-1})^T)$.

\medskip

{\bf Proposition 4 (\cite{3})\footnote{A.I. Sapelkin in \cite{5}
called this statement as Zeldich lemma. M.V. Zeldich in \cite{6}
called it  ``important and surprising'',
with what authors quite agree, in spite of brevity of the proof}}
{\it If there are no parallel pathes in $Q(S)$ then the forms
$\mathcal{T}_S$ and $f_S$ are equivalent over $\mathbb{Z}$.}

\medskip

Really, $\hat Q^{-1}\mathcal{A}_T(\hat Q^{-1})^T=\frac 12 \hat Q^{-1}(\hat
Q+\hat Q^T) (\hat Q^{-1})^T=\frac 12 [(\hat Q^{-1})^T +\hat
Q^{-1}]=A$.\hfill\rule{2mm}{2mm}

\medskip

Propositions~1, 2, 3, 4 imply

{\bf Corollary 3.} {\it Let $\Gamma(S)$ be acyclical and at least one of
the forms $f_S$, ${\mathcal T}_S$ be positive definite (it
holds if $\Gamma(S)$ is a Dynkin graph, see s.~4). Then the other form is
positive definite as well, and the
following statements are equivalent:

a) $S$ is antimonotonous;

b) $S$ is faithful;

c) ${\rm St}(S)\neq \varnothing$.}

\medskip

We denote $I(s_i)$  for $s_i\in S$ the number of edges of graph
$\Gamma(S)$, having $s_i$ its terminus.

$s_i$ is {\it a terminal point} if $I(s_i)\leq 1$, $s_i$ is
{\it a branch point} if $I(s_i)\geq 3$. $s_i$ is {\it a junction point},
if it is either a terminus of at least two arrows, or an origin of at least
two arrows of quiver $Q(S)$. We denote $S^\times$ the set of junction points.

\medskip

{\bf Example 4.}

\vspace{2mm}

$$
S=\mbox{\raisebox{-10mm}[0pt][0pt]{\includegraphics{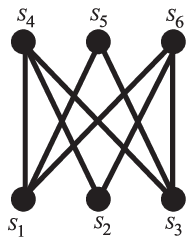}}}
$$

\vspace{7mm}

Here $S^\times =S$, $C(S)\cup \, \widetilde{\rm St}(S)=\varnothing$.
Really, values of $f_S$ can be negative ($f_S(1,1,1,-1$, $-1,-1)=-2$),
consequently, by Proposition~3,
$\widetilde{\rm St}(S)=\varnothing$, but ${\rm St}(S)\ni (1,2,1,1,2,1)$,
$|A|=-48$.
Proposition~2 implies $C(S)=\varnothing$.

We will assume further (except Appendix) graph $\Gamma(S)$ to be connected.

We remind that a cycle in graph $\Gamma$ is a sequence $\{s_1,\dots s_m\}$
of different vertices $s_i$ of graph $\Gamma$ where $m\geq 3$ such that a
vertex $s_i$ is connected with $s_{i+1}$ for $i=\overline{1,m}$, and
a vertex $s_m$ is connected with $s_1$. We call a cycle $\{s_1,\dots,s_m\}$
{\it simple}, if there are no other vertices between $s_1,\dots,s_m$.
We call graph $\Gamma$ and poset $S$ such that $\Gamma(S)=\Gamma$  {\it
cyclical}, if $\Gamma$ contains a cycle, and {\it acyclical} in opposite
case. it is easy to see that a cyclical graph contains a simple cycle, and,
correspondingly, a cyclical poset $S$ contains a subset $S'$ such that
$\Gamma(S')=\widetilde{A}_m$. It is clear that if $S$ is acyclical then
$Q(S)$ has no parallel pathes.

If $\Gamma(S)=\Gamma(\overrightarrow S)$ then $Q(\overrightarrow S)$ can
be obtained from $Q(S)$ by ``reorientation''
(i.~e. by alternation of direction)
of several arrows. If $\Gamma(S)$
is acyclic and quiver $\overrightarrow Q$ is obtained by reorientation of
arrows from $Q(S)$, then there exists $\overrightarrow S$ such that
$\overrightarrow Q=Q(\overrightarrow S)$.

We call quiver $Q=Q(S)$ and poset $S$
{\it standard}, if
$I(s_i)=2$ implies that $s_i$ is origin of one arrow and terminus also
of one arrow, and $I(s_i)\neq 2$ implies that $s_i$ is either origin of
$I(s_i)$ arrows or terminus of $I(s_i)$ arrows ($i=\overline{1,n}$).
It is easy to see that exactly one
standard poset is attached to each acyclical graph
up to antiisomorphism. If $S^*$ is antiisomorphic to $S$ then
$\Gamma(S)=\Gamma(S^*)$, and $Q(S^*)$ is obtained from $Q(S)$
by reorientation of all arrows.

If $\varphi$ is an arrow of $Q(S)$,
then we denote $S(\varphi)$ poset, obtaining from $S$
by overturn of an arrow $\varphi$.
$A_\varphi$ is matrix of $f_{S(\varphi)}$.
It is clear that $\mathcal{A}(S(\varphi))=\mathcal{A}(S)$.

We call a point $s_m\in S$ {\it Dynkin}, if there exists $m$-Dynkin vector
for form $f_S$.

\medskip

{\bf Remark 2.} Function $f_S$ is $m-$isolated in sense of s.1, if
$s_m$ is not comparable with other points of $S$. Therefore for connected
$S$ the condition of non-isolateness of $f_S$ in Lemma~3 holds
automatically.

\medskip

{\bf Lemma 6.}
{\it Let $\Gamma(S)$ be acyclical, $s_i\overset{\varphi}\rightarrow
s_j\in Q(S)$; $d\neq 0$ is such a vector that $d'_i=2(dA)_i=0$,
$d'_j=2(dA)_j=0$. Then there exists a vector $\widehat{d}\neq 0$ such that
$dA=\widehat{d}A_\varphi$.}

\medskip

It follows from the proof of Proposition~4 that
$\widetilde{Q}^{-1}\mathcal{A}(\widetilde{Q}^{-1})^T=A$,
$\widetilde{Q}^{-1}_{\varphi} \mathcal{A}_{\varphi}
(\widetilde{Q}^{-1}_{\varphi})^T = A_{\varphi}$,
$\widetilde{Q}_\varphi^{-1}\widetilde{Q}A
\widetilde{Q}^T(\widetilde{Q}_\varphi^T)^{-1}
=A_\varphi$ ($\mathcal{A}_{\varphi}=\mathcal{A}$).
Put
$\widehat{d}=d\widetilde{Q}^{-1}\widetilde{Q}_\varphi$, $(\widehat{d}\neq 0)$,
$\widehat{d}A_\varphi=d\widetilde{Q}^{-1}\widetilde{Q}_\varphi
\widetilde{Q}_\varphi^{-1}\widetilde{Q}A\widetilde{Q}^T(Q_\varphi^T)^{-1}
= dA\widetilde{Q}^T(\widetilde{Q}_\varphi^T)^{-1}$.

Put (in view of $d'_i=d'_j=0$),
$dA=\sum\limits_{k\not\in\{i,j\}}\alpha_ks_k=b$. We shall prove that
$b\widetilde{Q}^T(\widetilde{Q}_\varphi^T)^{-1}=b$, i.~e. that
$s_k\widetilde{Q}^T(\widetilde{Q}_\varphi^T)^{-1}=s_k$, i.~e. that
$s_k\widetilde{Q}^T=s_k\widetilde{Q}_\varphi^T$, that follows from
the definition of $\widetilde{Q}$ and from $k\not\in\{i,j\}$.\hfill
\rule{2mm}{2mm}

\medskip

{\bf Lemma  7.} {\it If $\Gamma(S)$ is acyclical, $s_t$ is a Dynkin
terminal point of $S$ then it is Dynkin point for poset
$\overrightarrow S$ if $\Gamma(\overrightarrow S)=\Gamma(S)$.}

\medskip

If $\overrightarrow{S}=S(\varphi)$
then the statement follows from Lemma~6.
Considering the general case ($\overrightarrow S$ is not
$S(\varphi)$), we remark firstly that if $S^*$ and $S$ are antiisomorphic
then $f_S=f_{S^*}$ and $s_t$ is a Dynkin point also for $S^*$.

We denote $\psi$ (resp. $\hat\psi$)
unique arrow $Q(S)$ (resp. $Q(\overrightarrow{S}$))
for which $s_t$ is either terminus or origin.
Then the condition $s_t\not\in\{s_i,s_j\}$
is equal to $\varphi\not=\psi$.

Without loss of generality suppose that $\psi$
in $Q(S)$ and $\hat\psi$ in $Q(\overrightarrow S)$
have the same orientation, otherwise pass to $\overrightarrow S^*$. In this
case we can pass from $S$ to $\overrightarrow S$, returning
several arrows, different from $\psi$, and, therefore, partial cases
$\overrightarrow S=S(\varphi)$ (Lemma~6) and $\overrightarrow S=S^*$
considered by us imply the statement of lemma.\hfill\rule{2mm}{2mm}

{\bf 4.} Let $\Gamma$ be a connected acyclical graph
with one branch point and three terminal points.
$\Gamma$ is a union of three chains
$A_{n_1}$, $A_{n_2}$, $A_{n_3}$,
intersecting in a branch point $s_1$.
$\Gamma=A_{n_1}\cup A_{n_2}\cup A_{n_3}$,
$A_{n_1}\cap A_{n_2}=A_{n_1}\cap A_{n_3}=A_{n_2}\cap A_{n_3}=\{s_1\}$,
$|A_{n_j}|=n_j$, $j=\overline{1,3}$, $|\Gamma|=n_1+n_2+n_3-2$. We will
denote $\Gamma$ by $\Gamma( n_1,n_2,n_3)$ (graph does not change
if we permutate $n_j$).

All Dynkin graphs besides $A_n$ (i.~e. $D_n$, $E_6$, $E_7$, $E_8$)
and  extended Dynkin graphs $\widetilde{E}_6$, $\widetilde{E}_7$,
$\widetilde{E}_8$ have the form  $\Gamma(n_1,n_2,n_3)$. It is well-known
that $\Gamma(n_1,n_2,n_3)$ is a Dynkin graph if and only if
$n_1^{-1}+n_2^{-1}+n_3^{-1}>1$, and is an extended Dynkin graph when
$n_1^{-1}+n_2^{-1}+n_3^{-1}=1$.

Namely, $\Gamma(n_1,n_2,n_3)$ is $E_6$, $E_7$, $E_8$, $D_n$,
if  respectively $(n_1,n_2,n_3)=(3,3,2)$;    $(2,4,3)$; $(2,3,5)$;
$(1,1,n-2)$.
$\Gamma( m_1,m_2,m_3)$
is respectively $\widetilde{E}_6$, $\widetilde{E}_7$, $\widetilde{E}_8$.
for $(m_1,m_2,m_3)=(3,3,3)$; $(2,4,4)$; $(2,3,6)$.

We fixed here the numeration of $m_j$ and $n_j$
so that $m_1\leq  m_2\leq m_3$
and for $E_n$ $(n=6,7,8)$,
$n_1=m_1$, $n_2=m_2$, $n_3=m_3-1$.

Remark, that (in all cases) $m_3$ is divided by $m_1$ and $m_2$.

\medskip

{\bf Proposition 5.} {\it If $\Gamma(S)=\Gamma(n_1,n_2,n_3)$
is Dynkin graph or extended Dynkin graph (i.~e. $D_n$, $E_6$,
$E_7$, $E_8$, $\widetilde{E}_6$, $\widetilde{E}_7$ or $\widetilde{E}_8$)
then $S$ contains terminal Dynkin point.}

\medskip

By Lemma~7 we can, without loss of generality suppose $S$
to be standard.

For any $\Gamma(m_1,m_2,m_3)$ $(I(s_1)=3)$ we construct a vector $\widetilde d$,
putting $\widetilde{d}_1=-m_3$; $\widetilde{d}_i=\frac{m_3}{m_j}$,
for $s_i\in A_{n_j}$, $i\not=1$.
It is easy to see that $\widetilde{d}'_i=0$, for $i\not=1$,
$\widetilde{d}'_1=m_3(1-m_1^{-1}-m_2^{-1}-m_3^{-1})$.
If, moreover,  $\Gamma(S)$ is extended Dynkin graph then $\widetilde{d}\in
\mathbb{Z}_n$ and $\widetilde{d}'_1=0$, i.~e. $\widetilde d$ is
in this case a $i$-Dynkin vector for any $i$.

Let $\Gamma(S)$ be $E_n$, $|S|=n$.
$\widetilde S$ is such standard poset that $\Gamma(\widetilde
S)=\widetilde{E}_n$, $|\widetilde S|=n+1$, $S\subset \widetilde S$,
$\widetilde S\setminus S=\{s_{n+1}\}\subset A_{n_3}$. We construct for $S$
a Dynkin vector $d$, modifying Dynkin vector $\widetilde d$ for
$\Gamma(\widetilde S)$. Put $d_i=\widetilde{d}_i$ for $i<n$, and $d_n=2$
$(=\widetilde{d}_n+\widetilde{d}_{n+1})$, $d'_n=1$; ($d'_i=0,\;
i=\overline{1,n-1}$).

Let
\[
\Gamma(S)=D_n=\raisebox{-9mm}[+10mm][0pt]{\mbox{\includegraphics{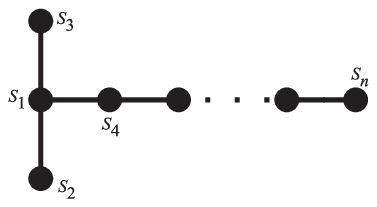}}},
\]

\vspace{5mm}

\noindent
then $w=(w_1,\ldots,w_n)$
where $w_1=-2$, $w_2=w_3=1$, $w_n=2$, $w_i=0$, $i\not\in    \{1,2,3,n\}$
is $n$-Dynkin vector $s_n$ $(w'_n=2)$.\hfill \rule{2mm}{2mm}

We write out evidently Dynkin vector for standard $S$
if $\Gamma(S)=E_6,E_7,E_8$
$$
\mbox{\includegraphics{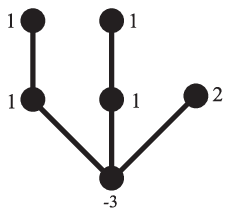}}\qquad\mbox{\includegraphics{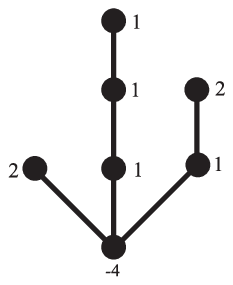}}\qquad
\mbox{\includegraphics{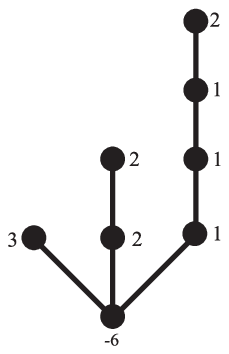}}
$$

{\bf Example 5.} Dynkin vector for standard poset $S$
such that $\Gamma(S)=\widetilde{D}_n$, $n>4$
(for $n=4$ see Example~2)
has the following form (all points are Dynkin points)

$$
\mbox{\includegraphics{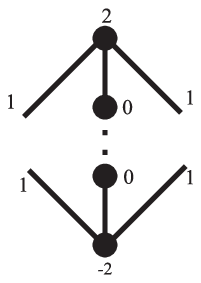}}
\qquad 
$$

{\bf 5.} Let poset $V=\raisebox{-9mm}[+10mm][10mm]{\mbox{\includegraphics{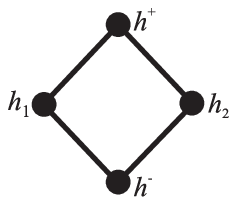}}}$,
poset $W^{2k}=\{s_1^-,\ldots,s_k^-,s_1^+,\ldots,s_k^+\,|\, s_i^-<s_i^+, s_i^-<s_{i+1}^-,
s_k^-<s_1^+, i=\overline{1,k}\}$, $k>1$,
in particular,
$W^4=\raisebox{-9mm}[+10mm][10mm]{\mbox{\includegraphics{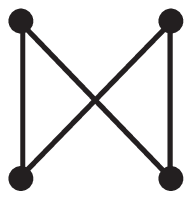}}}$
(see Example~3).

\medskip

{\bf Lemma 8.} {\it  If poset $S$ is cyclical and each $S'\subset S$
is acyclical then $S$ is either $V$, or $W^{2k}$ $(k\geq 2)$.}

\medskip

Without loss of generality suppose that $\Gamma(S)$
is a simple cycle $\widetilde{A}_n=
\raisebox{-3mm}[7mm][4mm]{\mbox{\includegraphics{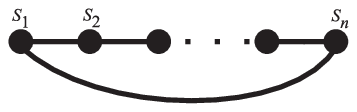}}}$.
Let
$S\setminus S^{\times}\ni s$, then $s^-<s<s^+$ and in $\Gamma(S)$ $s^-$ and
$s^+$ are connected with $s$ by edges. Then consider $S'=S\setminus\{s\}$.
If there is an edge $s^-$---$s^+$ in graph $\Gamma(S')$ then $\Gamma(S')$
is a cycle (that contradicts the condition of lemma).

If $s^-$ and $s^+$ are not connected by edge in $\Gamma(S')$
then there is a point $\bar   s\not=s$ in $S$ such that
$s^-<\bar s<s^+$, $\bar s>\!\!\!\!< s$ (i.~e.
$\bar s$ and $s$ are not comparable) because in the opposite case
$s^-$ and $s^+$ should be not connected with $s$ in $\Gamma(S)$.
And so $\{s^-,s^+,s,\bar s\}=V$ and $S=V$.

So, if $S\neq V$ then $S^{\times}=S$ and then
(since $\Gamma(S)=\widetilde{A}_n$) it is easy to see that
$S=W^{2k}$ for some $k>1$.\hfill \rule{2mm}{2mm}

\medskip

{\bf Lemma 9.} {\it If $S\supseteq V$ and $S\not\supset W^4$ then
$C(S)\not=\varnothing$.}

\medskip

An arbitrary vector $v\in {\mathbb R}_n$ can be considered as a function
on $S$ with values in ${\mathbb R}$.

Let $v: S\to {\mathbb R}$. $v(h^{-})=v(h^+)=-1$,
$v(h_1)=v(h_2)=1$, $v(t)=0$ for $t\in S\setminus V$; $v\in H_n$.
We prove that $v'(s)\leq 0$ for $s\in S$.
$v'(h^-)=v'(h^+)=-1$,
$v'(h_1)=v'(h_2)=0$.
If $t$ is not comparable niether with $h_1$ nor with $h_2$
then it is clear that $v'(t)\leq 0$.
If $t$ is comparable only with one of $h_1$, $h_2$,
then it is comparable either with $h^-$ or with $h^+$,
and also $v'(t)\leq 0$.
Let $t$ be comparable with $h_1$ and with $h_2$.
Suppose $t<h_1$ (a case $t>h_1$ is analogous) then $t<h_2$.
($h_1<t<h_2$ is impossible, so $t<h^+$).
Then if $t$ is comparable also with $h^-$ then $v'(t)=0$,
and in the opposite case $S\supset
W^4=\{t,h_2,h^-,h_1\}$.~\hfill~\rule{2mm}{2mm}

\medskip

{\bf Lemma 10.} {\it If $S\supseteq W^{2k}$ $(k\geq 2)$
and form $f_S$ is positive semidefinite then $C(S)\not=\varnothing$.}

\medskip

Let $t\in T=S\setminus W^{2k}$.
Put $S^-(t)=|\{s_i^-\,|\, t<s_i^-\}|\cup \{s_i^-\,|\, t>s_i^-\}|$,
$S^+(t)=|\{s_i^+\,|\, t<s_i^+\}\cup \{s_i^+\,|\, t>s_i^+\}|$.
We prove that positive semidefiniteness of $f_S$ implies
$S^-(t)=S^+(t)$.

Really, let $S^-(t_0)>S^+(t_0)$ for fixed $t_0\in T$
(a case $S^-(t_0)<S^+(t_0)$ is analogous).
Consider $x: S\to \mathbb{R}_n$ where $x(s_i^-)=-1$, $x(s_i^+)=1$,
($i=\overline{1,k}$), $x(t_0)=\varepsilon$, $0<\varepsilon<1$
and $x(t)=0$ for $t\in T\setminus \{t_0\}$.
It is easy to see that $f_S(x)<0$.

Now consider vector $v: S\to {\mathbb R}_n$
where $v(s_i^-)=-1$, $v(s_i^+)=1$,
$v(t)=0$ for $t\in T$.
$S^{-1}(t)=S^+(t)$ $(t\in T)$ implies
$v'(s)=0$ for any $s\in S$. It is clear, that $v\in H_n$ and
therefore $v\in C(S)$\hfill \rule{2mm}{2mm}

\medskip

{\bf Proposition 6.} {\it  If $S$ is an antimonotonous poset and form $f_S$
is positive semidefinite then $\Gamma(S)=A_n$.}

\medskip

If $S$ is cyclical, then Lemmas 8,9,10 imply the statement. If $S$ is
acyclical, then by Proposition~4 the Tits form $\mathcal{T}_S$ is
positive semidefinite, so $\Gamma(S)$ is one of $A_n$, $D_n$, $E_6$, $E_7$,
$E_8$, $\widetilde{D}_n$, $\widetilde{E}_6$, $\widetilde{E}_7$,
$\widetilde{E}_8$ ($\Gamma(\widetilde{A}_n)$ is cyclical). If
$\Gamma(S)\neq A_n$ then Proposition~5, Examples~2 and 5 and Lemma~7 imply
the existence of the Dynkin point and, by Lemma~3 (and Remark~2), $C(S)\neq
0$.

{\bf 6.} We consider now $\Gamma(S)=A_n$. In this case poset $S$ up
to antiisomorphism is defined by its order and
by subset $S^\times$ consisting of junction points (see s~.3). It is clear
that $S^\times=\emptyset$ if and only if $S$ is a chain.

Let poset $W^{k,k+1}=\{s_1^-,\ldots,s_k^-,s_1^+,\ldots,s_{k+1}^+\,
|\, s_i^-<s_i^+, s_i^-<s_{i+1}^+,i=\overline{1,k}\}$,
$W^{k+1,k}=\{s_1^-,\ldots,s_{k+1}^-,s_1^+,\ldots,s_{k}^+\,
|\, s_i^+>s_i^-, s_i^+>s_{i+1}^-,i=\overline{1,k}\}$.

\medskip

{\bf Lemma 11.} {\it If $\Gamma(S)=A_n$ and
$\Gamma(S)\supseteq W$ of the form $W^{k,k+1}$ (resp. $W^{k+1,k}$),
moreover $s_1^+,s_{k+1}^+\not\in S^\times$ (resp.
$s_1^-,s_{k+1}^-\not\in S^\times$) then $C(S)\not=\varnothing$.}

\medskip

Let, for the sake of definiteness, $S\supset W^{k,k+1}$. Consider vector
$v$ such that $v(s_i^-)=-2$, $i=\overline{1,k}$;
$v(s_i^+)=+2$, $i=\overline{2,k}$;
$v(s_1^+)=v(s_{k+1}^+)=1$;
$v(t)=0$ for $t\in S\setminus W^{k,k+1}$. We prove that $v\in C(f)$.

Indeed, $v\in H_n$; $v'(s_1^-)=v'(s_k^-)=-1$; $v'(s_i^-)=0$ for
$i=\overline{2,k-1}$; $v'(s_i^+)=0$ for $i=\overline{1,k+1}$.
Absence of branch points implies that if $t\not\in W$ is comparable with
$w\in W$ then $w\in \{s_1^+,s_{k+1}^+\}$. If $t$ is comparable with both
$s_1^+$, $s_{k+1}^+$ then $S$ would be cyclical. $s_1^+$, $s_{k+1}^+$
$\not\in S^{\times}$ implies that it can be only $t > w$. Therefore
every $t$ either is comparable exactly with one $s_i^-$ and one
$s_i^+$, or is not comparable with any $w\in W$, hereof $v'(t)\leq 0$.
\hfill \rule{2mm}{2mm}

We call a poset $\zeta$ {\it a wattle} \cite{2} if it is the union of not
intersecting chains $Z_i$, $|Z_i|\geq 2$, $i=\overline{1,t}$, $t>1$, in
which the minimal element of $Z_i$ is less than the maximal element of
$Z_{i+1}$ and there are no other comparisons between elements of different
$Z_i$. $\Gamma(\zeta)=A_n$. According to \cite{2} we denote $\zeta=\langle
n_1,\ldots,n_t\rangle$ where $n_i=|Z_i|$.

For a poset $S$ we consider a disconnected subgraph
$\Gamma(S^{\times})$ of $\Gamma(S)$. Denote $S_i^{\times}$ its connected
components.

\medskip

{\bf Lemma 12.}
{\it A poset $S$ where $\Gamma(S)=A_n$ is either a chain or a wattle
if (and only if) the orders of all $S_i^\times$ are even.}

\medskip

If $S$ is a wattle then the statement is evident (and we will not use it).
We will prove the converse statement induction by $|S|$. The base is
evident. Let $|S|=n+1$. $\Gamma(S)=\cdots s_{n-1}$---$s_n$---$s_{n+1}$
where $s_{n+1}$ is a terminal point (therefore $s_{n+1}\not\in
S^\times$). For the sake of definiteness we will suppose that
$s_n>s_{n+1}$. So $s_{n+1}$ is minimal. Put $S'=S\setminus\{s_{n+1}\}$
and $S''=S\setminus\{s_{n+1},s_{n}\}$. We have two
possibilities: 1)$s_{n-1}>s_n$; 2)$s_{n-1}<s_n$.

1) $s_n\not\in S^\times$, $(S')^\times=S^\times$. By assumption of
the induction $S'$ is a wattle in which $s_n$ is a minimal
terminal point. It is clear that $S$ is either a wattle or a chain. (If
$S'=\langle n_1,\dots,n_t\rangle $ then $S'=\langle
n_1,\dots,n_t+1\rangle$).

2) $s_n\in S^\times$ ($S'$ does not satisfy the induction's
assumption!). $s_n\in S^\times_p$, $|S_p^\times|\equiv 0 (\mod 2)$.
So $s_{n-1}\in S_p^\times\subset S^\times$. $s_{n-1}$ is a terminal point of
$S''$ and so $s_{n-1}\not\in(S'')^\times$. If
$S^\times=\bigcup\limits_{i=1}^{p}S_i^\times$ then
$(S'')^\times=\bigcup\limits_{i=1}^{p-1}S_i^\times
\cup(S^\times_p\setminus\{s_n,s_{n+1}\})$.

Consequently $S''$ satisfies the induction's assumption and
so it is either a chain or a wattle, in which $s_{n-1}$ is a minimal
point. If $S''=\langle n_1,\dots n_t\rangle $ then
$S=\langle n_1,\dots,n_t,2\rangle$.\hfill\rule{2mm}{2mm}

\medskip

{\bf Proposition 7.} {\it If a form $f_S$ is positive semidefinite and
$C(S)=\varnothing$, then $S$ is either a chain or a wattle.}

\medskip

Proposition~6 implies $\Gamma(S)=A_n$.

If $S$ is neither a chain nor a wattle then by Lemma~12 there exists
$S^\times_p$ such that $|S_p^\times|\equiv 1(\mod 2)$. It is easy to see
that $S_p^\times$ is $W^{k,k+1}$ or $W^{k+1,k}$. Let, for the sake of
definiteness $S_p^\times=W^{k,k+1}=
\raisebox{-8mm}[7mm][8mm]{\mbox{\includegraphics{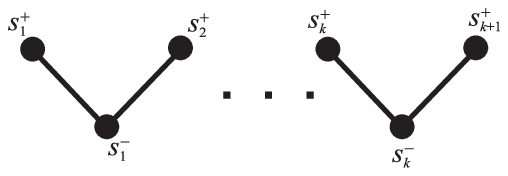}}}$.
The facts that $s_1^+\in S^\times$ and $S_p^\times$ is a connected component
in $S^\times$, imply that there exists $s_0^-\in S\setminus S^\times$
such that $s_0^-\rightarrow s_1^+$. Analogously there exists $s_{k+1}^-\in
S\setminus S^\times$, $s_{k+1}^-\rightarrow s_{k+1}^+$. It is easy to see
that $S_p^\times\cup\{s_0^-,s_{k+1}^-\}=W^{k+2,k+1}$ and
$C(S)\neq\emptyset$ by Lemma~11.\hfill\rule{2mm}{2mm}

\medskip

Example~4 (s.~3) implies that the condition of positive semidefiniteness of
$f_S$ can not, generally speaking, be excluded.

\medskip

{\bf Hypothesis.} {\it If $S$ is acyclical and $\Gamma(S)\neq A_n$ then
$C(S)\not=\varnothing$.}

\medskip

In some cases the existence of $v\in C(f)$ for acyclical $S$
is evident. We give however the example of an acyclic poset $S$ and
$v\in C(S)$ which we constructed only by means of computer.

\medskip

{\bf Example 6.}
$$
\raisebox{-65mm}[0mm][0mm]{\mbox{\includegraphics{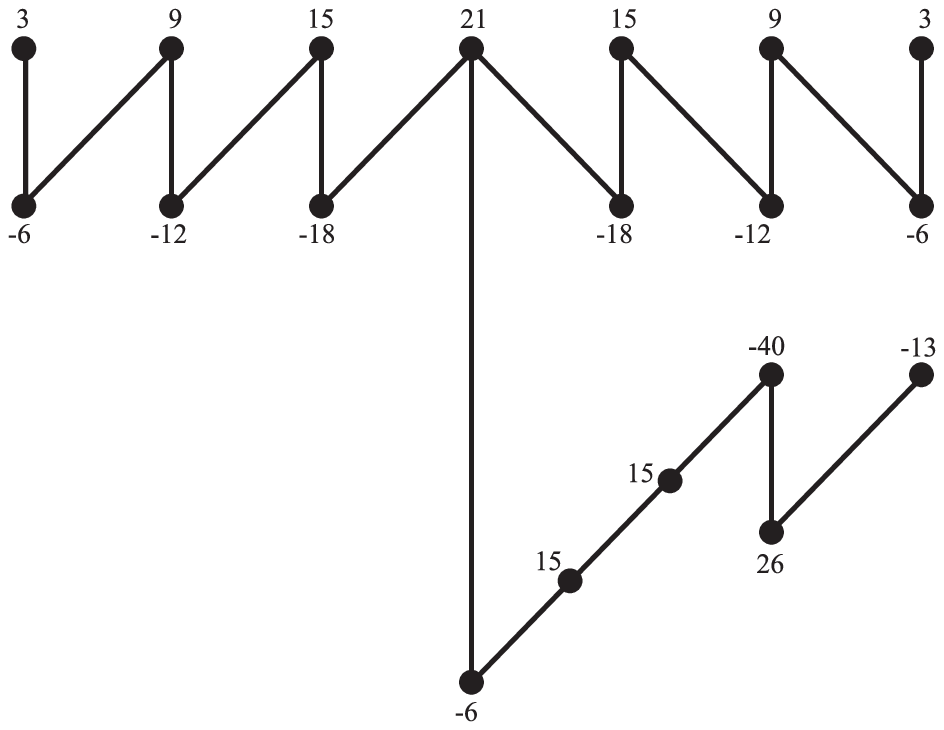}}}
$$

\vspace{65mm}

{\bf 7.} Let $\zeta=\langle n_1,\ldots,n_t\rangle$ $(t>1)$ be a wattle,
where $n_i=|Z_i|$, $\sum\limits_{i=1}^t n_i=n$, $n_i>1$, $i=\overline{1,t}$.

In \cite{2} the minimal points of chains $Z_i$, $i=\overline{1,t-1}$
are denoted by $z_i^-$, and maximal points of chains $Z_i$,
$i=\overline{2,t}$ are denoted by $z_i^+$, $z_i^-<z_{i+1}^+$. The rest
(i.~e. not junction) points are called {\it common} (including
the maximal point of chain $Z_1$ and a minimal point of chain $Z_t$). They
are compared only with points from its chain.

A  {\it width} $\omega(S)$ of a
partially ordered set $S$
is the maximal number of its pairwise uncomparable elements.
We attach to any poset $S$
a rational number $r(S)=\frac{n+1}{t}-1$
where $n=|S|$, $t=\omega(S)$.
If $S$ is a chain then $w(s)=1$, $r(S)=n$.

It is clear that there exist many wattles with the same $r$. We prove below,
however, that to any noninteger $r>1$ it is corresponded exactly one
(uniform in sense~\cite{2}) $P$-faithful ($=$ antimonotonous, see
Corollary~3) wattle, which we will call here $r$-wattle and denote
$\zeta(r)$.

For a positive rational $a$ we put $\{a\}=a-[a]$. Let $r$ be a positive
noninteger rational number greater than 1. $q/t$ is the representation
of $\{r\}$ in the form of irreducible fraction. We write out the
sequence of integers $n_1,\ldots,n_t$, which
will be the orders of sets $Z_i$ in $\zeta(r)$. Put $n_1=n_t=[r]+1$,
$n_i=[ri]-[r(i-1)]+1$ for $i=\overline{2,t-1}$. It is clear, that
$[\{r\}i]-[\{r\}(i-1)]$ is either 1 or 0. Therefore $n_i$ is either $1+[r]$
or $2+[r]$. The number of those $i$ for which $n_i=2+[r]$ is $q-1$;
$n=t([r]+1)+q-1$, $r(\zeta(r))=r$

Remark, that $r$-wattles are uniform in sense of \cite{2}
and conversely.

Thus, we corresponded to each noninteger rational number
$r>1$ a wattle $\zeta(r)$. We can consider also integer numbers,
putting for natural $r$ that $\zeta(r)$ is a chain of length
$r$. We will call $r$-sets all posets of form $\zeta(r)$ $r\geq 1$
(i.~e. uniform wattles and chains).

\medskip

{\bf Theorem.} {\it Let form $f_S$ be
positive semidefinite ($\Gamma(S)$ be connected). Then $C(S)=\varnothing$,
if and only if $S$ is $r$-set.}

\medskip

If $r$ is integer then the statement is evident (see~\cite{2}).
Therefore as a matter of fact we need to prove, counting
Proposition~7, that $C(\zeta)=\varnothing$ if and only if $\zeta$
is $r$-wattle.

We attach to $r$-wattle $\zeta(r)$ a vector $x: \zeta \to {\mathbb R}^+$,
$x(s)=1$ for $s\in\zeta\setminus\zeta^\times$, $x(z_i^-)=\{ir\}$,
$x(z_i^+)=1-x(z_{i-1}^-)$.
$(q,t)=1$ implies $x(s)>0$ for any $s\in \zeta(r)$.

It is posiible to check that $x\in {\rm St}(\zeta)$
either spontaneously, using the definition of
${\rm St}(\zeta)$, or using the following lemma.

\medskip

{\bf Lemma 13 (see~\cite{2}, Lemma~5).} {\it  Vector $x: \zeta \to {\mathbb
R}$ is contained in ${\rm St}(\zeta)$ if and only if there exist such
positive $\alpha$, $\beta$ that

$1)$ $x(s)=\alpha$ for $s\in \zeta\setminus \zeta^\times$ (we can suppose
$\alpha=1$ multiplying $x$ by $\lambda\in{\mathbb R}^+$);

$2)$ $x(z_i^-)+x(z_{i+1}^+)=\alpha$ for $i=\overline{1,t-1}$;

$3)$ $\sum\limits_{s\in Z_i} x(s)=\beta$, $i=\overline{1,t}$.}

\medskip

Proof of the lemma is almost evident. Remark merely that at first
2) should be proved (it follows from $\frac{\partial f_S}{\partial
z_i^-}(x) = \frac{\partial f_S}{\partial z_i^+}(x)$, $i=\overline{2,t-1}$),
and then 1).

Vector $x$ constructed above
evidently satisfies the conditions 1), 2) of Lemma~13.It is easy to
check (for $\alpha=1$) that

\begin{gather}
\sum\limits_{s\in Z_i} x(s)=r (i=\overline{1,t});\label{4*}\\
x'(s)=1+r;\label{5*}\\
\sum\limits_{s\in \zeta(r)} x(s)=tr,\label{6*}
\end{gather}

Thus, $x\in{\rm St}(\zeta)$, so $\widetilde{\rm St}(\zeta)\neq \varnothing$,
and $C(S)=\varnothing$ by Corollary~3.

It remains to prove that any
$P$--faithful wattle $\zeta$ is
$r$-wattle (where $[z]=|Z_1|-1$, $\{r\}=x(z_1^-)$). This follows from the
next statement.

{\bf Lemma 14.}
{\it Let $\zeta=\langle z_1,\ldots,z_t\rangle$
and $\hat \zeta=\langle \hat z_1,\ldots,\hat z_t\rangle$ are two
$P$-faithful wattles, $x\in{\rm St}(\zeta)$, $\hat x\in {\rm St}(\hat
\zeta)$ $(\hat{\alpha}=\alpha=1)$. Then, if $Z_1=\hat Z_1$ and $x(s)=\hat
x(s)$ for $s\in Z_1=\hat Z_1$ then $\zeta=\hat \zeta$ and $x(s)=\hat x(s)$
for $s\in\zeta$.}

It is sufficient to prove that if
$m\leq\max\{t,\hat t\}$,
then $z_i=\hat z_i$ for $i\leq m$
and $x(s)=\hat x(s)$
for $s\in \cup_{i=1}^m Z_i$.
Lemma~13 implies this by induction with respect to $m$
(see~\cite{2}). \hfill\rule{2mm}{2mm} \hfill\rule{2mm}{2mm}

There was introduced the numerical
function $\rho(r)=1+\frac{r-1}{r+1}$ where $r\in
{\mathbb N}$ \cite{2}. We spread this definition on the case $r\geq
1$ is rational. Put $\rho(r_1,\ldots,r_t)=\sum\limits_{i=1}^t\rho (r_i)$.
If $Z_n$ is a chain of order $n$ then $P(Z_n)=\rho(n)$ \cite{2}. Let
$\zeta(r)$ be a wattle. Vector $\overline{x}=(tr)^{-1}x\overset{(6)}{\in}
P_n \cap {\rm St} (\zeta(r))$ (where $x$ is a vector
constructed in the proof of the theorem).

\begin{gather*}
P(\zeta(r))=f_{\zeta(r)}^{-1}(\overline{x})=
(tr)^2f_{\zeta(r)}^{-1}(x)\overset{(3)}{=}\\
=2(tr)^2 \left(\sum\limits_{s\in\zeta(r)}x'(s)x(s)\right)^{-1}
\overset{\eqref{5*},\; \eqref{6*}}{=}
\frac{2t^2r^2}{(1+r)tr}=\frac{2tr}{1+r}=t\rho(r).
\end{gather*}

This formula is true if $t=1$ (i.~e. in the case of chain). For
any positive rational $r=\frac{l}{t}$ $((l,t)=1)$
$t\rho(r)=\frac{2lt}{l+t}$. Introduce the function $P(r)=\frac{2lt}{l+t}$
(for $n\in {\mathbb N}$ $P(n)=\rho(n)$). Thus, for any $r\geq 1$

\begin{equation}
P(\zeta(r))=t\rho(r)=P(r)\label{z4}.
\end{equation}

{\bf Appendix.} We call a poset $S$ {\it connected} if graph $\Gamma(S)$
is connected. The theorem and Corollary~3 imply that a connected
poset $S$ is $P$-faithful if and only if it is $r$-set \cite{2,3,4,5,6}. On
the other hand, \cite{2,3,4,5,6} imply our theorem only if $f_S$ is
positive definite (but not positive semidefinite).

We remind about role of $P$-faithful
posets in representation theory. We will write
$S=S_1\bigsqcup S_2$, if $S=S_1 \cup S_2$, $S_1\cap
S_2=\varnothing$ and elements $S_1$ are not comparable with elements
$S_2$. $S= Z_1 \bigsqcup \cdots\bigsqcup
Z_p$ is {\it primitive}, if $Z_i$ are chains, $i=\overline{1,p}$. We denote
it $(n_1,\ldots,n_p)$ if $n_i=|Z_i|$. The characterization of
antimonotonous disconnected posets follows from the theorem and Lemma~2.

Any poset $S=\bigsqcup\limits_{i=1}^p S_i$ where $S_i$ are connected
components. According to Lemma~5, $P(S)=\sum\limits_{i=1}^p P(S_i)$, and
if $S$ is primitive then
$P(S)=\sum\limits_{i=1}^p \rho(n_i)= \rho(n_1,\ldots,n_t)$.

A role of quadratic forms in the theory of representations of quivers and
posets is well-known~\cite{11}.

The norm of a relation $\|S,\leq\|=\inf\limits_{u\in \overline{P}_n}
f_S(u)$ was introduced in~\cite{1} on base of form $f_S$. In view
of Lemma~5 it is naturally to consider instead of $\|S,\leq\|$
the function $P(S)=\|S,\leq\|^{-1}$ \cite{1}.

\medskip

{\bf Proposition 8.} {\it
$S$ has finite (respectively tame) type if and only if
$P(S)<4$ (respectively $P(S)=4$).}

\medskip

With this point of view the Kleiner's list of the critical posets~\cite{7}
is the list of {\it $P$-faithful} posets $S_i$, for which $P(S)=4$.

4 posets of Kleiner's list are primitive:
\[
\mbox{(I)}.\quad  (1,1,1,1), \quad (2,2,2), \quad (1,3,3),\quad  (1,2,5),
 \quad \mbox{and the fifth is} \quad (4)\bigsqcup K,
\]
where $K=\raisebox{-4mm}[7mm][5mm]{\mbox{\includegraphics{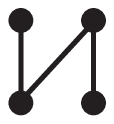}}}=
\langle 2;2\rangle=\zeta(1\frac 12)$.
It is easy to see that any chain is $P$-faithful, and s.~7 implies
that $K$ is also $P$-faithful $(P(K)=2,4)$. By Lemma~5 a
disconnected poset is $P$-faithful if and only if all its components
are $P$-faithful.

The list of critical sets \cite{9}:
\[
\mbox{(II)}.\quad (1,1,1,1,1),\quad (1,1,1,2),\quad (2,2,3),\quad
(1,3,4),\quad (1,2,6),\quad (6)\bigsqcup K
\]
\noindent can be characterized as the list of $S$, having the following
properties:

1) $P(S)>4$,

2) if $S'\subset S$ then $P(S')\leq 4$.

The following statement play central role in the theory of representations
of posets \cite{7}, \cite{9}.

{\it A poset $S$ is finitely represented (respectively
tame), if and only if $S$ does not contain subsets I (respectively
II).}

It has been naturally to suppose that all $P$-faithful poset are either
chains or some sets, for which $K$ is the least representative. This
was a reason to introduce $P-$faithful posets \cite{2}.

Now we show how the lists (I), (II) can be obtained using
characteristic of (connected) $P$-faithful sets and formula
\eqref{z4}. It is easy to check that $P(S)=4$ for $S\in \mbox{I}$ and
$P(S)>4$ for $S\in \mbox{II}$ (in view of Lemma~5 and formula \eqref{z4}).

We call a $P$-faithful poset $S$ {\it utmost} if $P(S)\geq
4$ and $P(S')\leq 4$ for any $S'\subset S$ (hereat $S'$ can be supposed
$P$-faithful).

\medskip

{\bf Lemma 15.} {\it Not primitive utmost $S$ is of the form $K\bigsqcup Z_m$
where $m$ is 4 or~5.}

\medskip

Let $S$ contains a connected component $\zeta(r)$ where
$\{r\}=\frac{q}{t}$, $t>1$, $q<t$, $(q,t)=1$. The characterization of
$P$-faithful posets implies $\omega(S)<4$ since in the opposite case
$S\supset S'=(2,1,1,1)$, $\rho(2,1,1,1)=4\frac{1}{3}>4$. Consequently
$t\leq 3$, moreover if $t=3$ then $\zeta(r)=S$.

Let $t=3$, $1\leq q\leq 2$. If $[r]\geq 2$ then $S\supset
S'=S\setminus\{z_1^-,z_3^+\}$. $S'$ is a primitive poset containing
$(2,3,2)$, $\rho(2,2,3)=4\frac{1}{6}$,
$\rho(S')>4$. If $[r]=1$ then either $r=1\frac{1}{3}$ or
$r=1\frac{2}{3}$. We obtain $\rho(r)\leq 1\frac{1}{4}$, then
$P(S)=3\rho(r)<4$ (see~\eqref{z4}).

Let $t=2$ and $S\neq
\zeta(1\frac{1}{2})\bigsqcup\hat{S}$. If $S=\zeta(r)$
then $P(S)=2\rho(r)<4$ since $\rho(r)<2$ for any $r$. So,
$S=\zeta(r)\bigsqcup\hat{S}$, $r>1\frac{1}{2}$, i.~e.
$r\geq 2\frac{1}{2}$, $\zeta(r)\supset\{\zeta(r)\setminus
z_1^-\}\supseteq(2,3)$. For $|\hat{S}|>1$, $S\supset S'$,
contains $(2,2,3)$ or $(1,1,2,3)$ and $P(S')>4$. Hence $|\hat{S}|=1$.
Then $[r]<3$, since otherwise $\zeta(r)\supset S'=\zeta(r)\setminus
z_1^-\supseteq(3,4)$ and $P(S'\bigsqcup(1))>4$ since
$\rho(3,4,1)>4$. For $[r]=2$ we obtain $P(S)<4$ since
$P(\zeta(2\frac{1}{2}))=2\cdot 1\frac{3}{7}$,
$P(\zeta(2\frac{1}{2})\bigsqcup(1))=2\frac{6}{7}+1$
(Lemma~5).

Let, finally, $S=\zeta(1\frac{1}{2})+\widehat{S}$
($S\neq\zeta(1\frac{1}{2})$) since $P(\zeta(1\frac{1}{2}))=2,4$),
then if $w(\widehat S)>1$ then $S\supset
S'=\{\zeta(1,\frac{1}{2})\setminus
z_1^-\}\bigsqcup(1,1)=(2,1,1,1)$, $\rho(S')>4$. If
$\widehat S=Z_m$ then for $m<4$ $\rho(m)<1,6$ and $P(S)<4$, and for $m>5$,
$S\supset S'=(\zeta(1\frac{1}{2})\bigsqcup Z_5)$,
$\rho(S')>4$, $P(K\bigsqcup Z_4)=4$,
$P(K\bigsqcup Z_5)=4\frac{1}{15}$.\hfill \rule{2mm}{2mm}

\medskip

{\bf Proposition 9.} {\it $P$-faithful $S$ is utmost if and only if $S\in
\mbox{\rm I}\cup \mbox{\rm II}$.}

\medskip

If $S$ is not primitive then Lemma~15 implies the statement.
Let $S$ is primitive. Then $w(S)>2$ and
$S\not\in\{(1,1,n),(1,2,2),(1,2,3),(1,2,4)\}$ (otherwise
$P(n_1,\dots,n_t)=\rho(n_1,\dots,n_t)<4$). In the rest cases
we can see that if $S\not\in \mbox{I}\cup \mbox{II}$ then
$S\supset S'\in \mbox{II}$, and if $S\in \mbox{I}\cup \mbox{II}$ then
$S\not\supset S'\in \mbox{II}$.\hfill\rule{2mm}{2mm}

\medskip

Propositions~8, 9 (in view of $P(S)=4$ for $S\in \mbox{I}$
and $P(S)>4$ for $S\in \mbox{II}$)
imply the main theorems \cite{7} and \cite{9}.

$P-$faithful posets, for which $P=4$, play important role in representation
theory. We don't know whether $P-$faithful posets with
$P=n>4$ play some analogous role. In~\cite{14} primitive posets with $P=5$
are written out. We give the example (probably unique) of a not primitive
poset $S=\zeta(3\frac{1}{2})\bigsqcup(17)$ for
which $P(S)=5$ (see \eqref{z4} and Lemma~3).

It is Example 4 (s.~3) where $C(S)=\varnothing$ but $S$ is not
$P-$faithful. We hope that studying of $C(S)$ can be interesting
for the representation theory. Remark,
that in~\cite{1} the norm $\|P\|$ of an arbitrary binary relation $P$ (on
finite set) and corresponding notion of $P$-faithfulness are defined (these
notions can be used for locally scalar representations (see~\cite{8}) in
Hilbert spaces). However in this case it is more complicate to review
$P$-faithful sets. Such investigation would be seemingly rather difficult
and interesting problem.

\end{document}